\documentclass[11pt]{article}
\usepackage{amsfonts, amsmath, amssymb, amsthm}

\usepackage[colorlinks=true,citecolor=blue,urlcolor=blue,linkcolor=blue,bookmarksopen=true]{hyperref}
\usepackage{aliascnt,cleveref}

\usepackage{algorithm}
\usepackage[noend]{algpseudocode}

\title{Twins in words and long common subsequences in permutations}
\author{Boris Bukh \and Lidong Zhou}

\makeatletter
\def\TODO{\@ifnextchar[{\TODO@with}{\marginnote{TODO}}}
\def\TODO@with[#1]{\marginnote{#1}}
\makeatother

\newtheorem{theorem}{Theorem}

\newaliascnt{question}{theorem}

\aliascntresetthe{question}
\crefname{question}{question}{questions}

\newaliascnt{lemma}{theorem}
\newtheorem{lemma}[lemma]{Lemma}
\aliascntresetthe{lemma}
\crefname{lemma}{lemma}{lemmas}

\newaliascnt{conjecture}{theorem}
\newtheorem{conjecture}[conjecture]{Conjecture}
\aliascntresetthe{conjecture}
\crefname{conjecture}{conjecture}{conjectures}

\theoremstyle{remark}
\newtheorem*{observation}{Observation}

\newcommand*{\abs}[1]{\lvert #1\rvert}                           
\newcommand*{\norm}[1]{\lVert #1\,\lVert}                        
\newcommand*{\eqdef}{\stackrel{\text{\tiny{def}}}{=}}            
\newcommand*{\E}{\mathbb{E}}                                     
\newcommand*{\Z}{\mathbb{Z}}                                     

\newcommand*{\W}{\mathcal{W}}                                    
\renewcommand*{\P}{\mathcal{P}}                                  
\newcommand*{\wordf}{\textsf}                                    
\newcommand*{\veps}{\varepsilon}                                 
\usepackage{stmaryrd}
\def\lbrac{\llbracket}
\def\rbrac{\rrbracket}
\newcommand*{\perm}[1]{\llbracket #1\rrbracket}

\newcommand*{\oline}[1]{\overline{#1\mathstrut}}                

\newcommand*{\alphab}[1]{\ifx 3#1
\{\wordf{0},\wordf{1},\wordf{2}\}
\else
\{\wordf{0},\wordf{1},\dotsc,\wordf{#1-1}\}
\fi}

\newcommand*{\alphabshort}[1]{\ifx 3#1
\{\wordf{0},\wordf{1},\wordf{2}\}
\else
\Sigma_k
\fi}

\renewcommand*{\alphab}[1]{[#1]}
\renewcommand*{\alphabshort}[1]{[#1]}

\DeclareMathOperator{\LCS}{LCS}                                  
\DeclareMathOperator{\LT}{LT}                                    
\def\LCSr{\LCS^{\operatorname{r}}}                               
\DeclareMathOperator{\rev}{rev}                                  
\def\emptyword{\emptyset}

\algloopdefx{IfOneline}[1]{\textbf{if} #1 \textbf{then}}         
\algloopdefx{RepeatTimes}[1]{\textbf{repeat} #1 times}           

\makeatletter
\newcommand*{\permt}[1]{\llbracket \def\l@amper##1{\ifx##1-{-}\else\ifx##1+\phantom{-}\else##1\fi\fi}\def\&{\quad\l@amper} \l@amper#1\rrbracket}
\makeatother

\date{}
\begin{document}
\maketitle
\begin{abstract}A large family of words must contain two words
that are similar. We investigate several problems where
the measure of similarity is the length of a common subsequence.

We construct a family of $n^{1/3}$ permutations on $n$ letters,
such that $\LCS$ of any two of them is only $cn^{1/3}$, improving
a construction of Beame, Blais, and Huynh-Ngoc. We relate
the problem of constructing many permutations with small $\LCS$
to the twin word problem of Axenovich, Person and Puzynina.
In particular, we show that every word of length $n$ over a $k$-letter
alphabet contains two disjoint equal subsequences of length $cnk^{-2/3}$.
Connections to other extremal questions on common subsequences are exhibited.

Many problems are left open.
\end{abstract} 
\section{Introduction}

This paper grew out of attempts to solve the twin word problem of Axenovich, Person
and Puzynina \cite{app_twins}. These attempts gave rise to several problems 
on common subsequences in words, some of which appear to us even more appealing.
However, we begin with the twin word problem, as it is the core that the other
problems link to.

A \emph{word} is a sequence of letters from some fixed finite alphabet.
Since the nature of the alphabet is not important to us, we will usually
use $\alphab{k}\eqdef \{1,2,\dotsc,k\}$ as a canonical $k$-letter alphabet. The
set of all words of length $n$ is thus denoted $\alphab{k}^n$. 
The $i$'th letter of a word $w$ is denoted $w[i]$.
A \emph{subword} of a word $w$ is a word consisting of several
consecutive letters from $w$. In contrast, a \emph{subsequence} in $w$
is a word consisting of letters from $w$ that are in order, but not
necessarily consecutive. For example, \textsf{radar} is a subsequence, but not a subword of 
\textsf{abracadabra}.

Two subsequences $w_1,w_2$ of a word $w$ are said to be \emph{twins}
if they are equal as words, and no letter from $w$ is in both $w_1$ and $w_2$. 
For example, the word \wordf{0110010010101101} contains twins of length $7$,
as seen from $\overline{\wordf{01}}\wordf{1}\underline{\wordf{0}}\overline{\wordf{0}}\underline{\wordf{10}}\overline{\wordf{0101}}\underline{\wordf{01}}\wordf{1}\underline{\wordf{01}}$,
where overlines and underlines indicate which of the two subsequences (if any) the letter is in.
Let $\LT(w)$ be the length of the longest twins contained in a word $w$. We also
define $\LT(k,n)\eqdef \min_w \LT(w)$, where the minimum is taken over all words $w\in \alphab{k}^n$. 
The twin problem of Axenovich--Person--Puzynina is to determine how large $\LT(k,n)$ is.
The following is a remarkable result:
\begin{theorem}[Axenovich--Person--Puzynina \protect\cite{app_twins}, Theorem 1]\label{thm:app_binary}
For words over a binary alphabet we have
\[
  \LT(2,n)\geq \tfrac{1}{2}n-o(n).
\]
\end{theorem}
Since the upper bound of $n/2$ is trivial, the result is tight. A very interesting
problem, which we do not address in this paper, is to determine the exact behavior 
of the $o(n)$ term. We refer the reader to \cite{app_twins} for the best known bounds. 
Instead we consider what happens as the alphabet size grows. First, we recall the results
of Axenovich--Person--Puzynina:
\begin{align}
      \LT(k,n)&\geq \tfrac{1}{k}n-o(n)&&\text{for all }k\geq 2,\label{eq:app_lb}\\
      \LT(5,n)&\leq 0.49n,\notag\\
      \LT(k,n)&\leq \left(\frac{e}{\sqrt{k}}-\frac{e^2}{k}+O(k^{-3/2})\right)n&&\text{for all large }k.
                          \label{eq:goestozero}
\end{align}
(The upper bound on $\LT(k,n)$ appearing in \cite{app_twins} is different; we have provided its asymptotic form.)

The upper bounds on $\LT$ from \cite{app_twins} use the union bound, which for $k=2,3,4$ yields only
the trivial bound $\LT(k,n)\leq n/2$. In particular, prior to the present work
the best bounds on $\LT(3,n)$ and $\LT(4,n)$ were $1/3\leq \lim_n \LT(3,n)/n\leq 1/2$,
and $1/4\leq \lim_n \LT(4,n)/n\leq 1/2$. We improve both lower and upper bounds. 

The first two results are improved lower bounds.
\begin{theorem}\label{thm:kthree}(Proof is in \cref{sec:regthree})
Let $k\geq 3$. Then we have $\LT(k,n)\geq \tfrac{1.02}{k}n-o(n)$. 
\end{theorem}
\begin{theorem}\label{thm:klarge}(Proof is in \cref{sec:commonseq})
For
the twin problem over a $k$-letter alphabet we have
\[ 
  \LT(k,n)\geq 3^{-4/3} k^{-2/3} n-3^{-1/3}k^{1/3}.
\]
\end{theorem}
While \cref{thm:klarge} is stronger than \cref{thm:kthree} for large $k$, together
they constitute an improvement over \eqref{eq:app_lb} for all $k\geq 3$. 

The next result is a new upper bound for $\LT(k,n)$. It is a tiny improvement over \eqref{eq:goestozero} for large $k$, 
but importantly it shows that $\LT(4,n)$ is not asymptotic to $n/2+o(n)$. 
(The question of whether or not $\LT(3,n)$ is asymptotic to $n/2+o(n)$ remains open.)
\begin{theorem}\label{thm:random_bound}(Proof is in \cref{sec:random_bound})
Suppose $k\geq 2$, and suppose $\alpha\in [1/k,1/2]$ is a constant such that the real number
\[
 (1-2\alpha)\log\frac{1}{1-2\alpha}-\alpha\log(\alpha^2k)-2\alpha\log \left(\frac{2}{1+\sqrt{1-1/k}}\right)+(1-2\alpha)\log(1-1/k)
\]
is negative. Then $\LT(k,n)\leq \alpha n$ for all sufficiently large $n$. In particular,
$\LT(4,n)\leq 0.4932n$, $\LT(5,n)\leq 0.48n$, and $\LT(k,n)\leq \left(\frac{e}{\sqrt{k}}-\frac{e^2+1/2}{k}+O(k^{-3/2})\right)n+o(n)$ for all large $k$.
\end{theorem}

In the next section we give a short proof of \cref{thm:klarge}.
The proof will illustrate a connection with several extremal 
problems on the longest common subsequences, which we discuss there as well.

\section{Common subsequences}\label{sec:commonseq}
A \emph{common subsequence} of two words $w$ and $w'$ is a word that is a 
subsequence of both $w$ and $w'$. A common subsequence of more than two words is defined
analogously. We denote the length of the longest common subsequence of a set $\W$
of words by $\LCS(\W)$. For notational sanity we write $\LCS(w,w')$ in place of $\LCS(\{w,w'\})$.

A special, but important class of words are permutations. A \emph{permutation} 
is a word in which each letter of the alphabet appears exactly once. We denote 
by $\P_k$ the set of all permutations on $k$ letters. Note that the group structure 
of permutations is of no concern in this paper; permutations are just words.
A slight generalization of permutations are multipermutations. Given a vector 
$\vec{s}=(s_1,\dotsc,s_k)\in \Z_+^k$, an \emph{$\vec{s}$-multipermutation} is a word in which
the letter $l\in \alphab{k}$ appears exactly $s_l$ times. Length of an $\vec{s}$-multipermutation
is $\norm{\vec{s}}\eqdef \sum_l s_l$. We denote the set
of all $\vec{s}$-multipermutations by $\P_{\vec{s}}$.

The basis for our proof of \cref{thm:klarge} is a result of
Beame and Huynh-Ngoc \cite[Lemma 5.9]{bh} which asserts that 
if $\pi_1,\pi_2,\pi_3\in \P_k$ are three permutations, then
for some pair $i<j$ we have $\LCS(\pi_i,\pi_j)\geq k^{1/3}$.
We require a slight generalization of that result:
\begin{lemma}\label{lem:multiperm_simple}
Let $\pi_1,\pi_2,\pi_3\in \P_{\vec{s}}$ be three $\vec{s}$-multipermutations of
length $\norm{\vec{s}}=n$. Then there is a pair $i<j$ such that
$\LCS(\pi_i,\pi_j)\geq n^{1/3}$. 
\end{lemma}
\begin{proof}
For each of $\pi_1,\pi_2,\pi_3$, replace the $j$'th occurrence of letter $l$ by the
letter $l_j$. An example: $\wordf{212331}$ becomes
$\wordf{2}_1\wordf{1}_1\wordf{2}_2\wordf{3}_1\wordf{3}_2\wordf{1}_2$. This
turns $\pi_i$ into a permutation $\pi_i'$ over an alphabet
of size $n$. Since $\LCS(\pi_i,\pi_j)\geq \LCS(\pi_i',\pi_j')$, the lemma follows from 
the result of Beame--Huynh-Ngoc stated above.
\end{proof}
\begin{proof}[Proof of \cref{thm:klarge}]
It suffices to show that every word of length $3k$ contains twins of length $(k/3)^{1/3}$, for
then the general result would follow by partitioning a word of length $n$ into intervals of length $3k$.
So assume that $w$ is of length $3k$.

For each $l\in \alphab{k}$ let $C_l$ be the number of copies of the letter $l$ in $w$.
Define a vector $\vec{s}$ by $s_l=\lfloor C_l/3 \rfloor$. Then $w$ contains three disjoint 
subsequences $w_1,w_2,w_3$, each of which is a $\vec{s}$-multipermutation.
Thus, by the preceding lemma there are $i<j$ such that 
\[
  \LCS(w_i,w_j)\geq \Bigl(\sum_{l\in \alphabshort{k}} \lfloor C_l/3\rfloor \Bigr)^{1/3}\geq \Bigl(\sum_{l\in \alphabshort{k}} (C_l-2)/3 \Bigr)^{1/3}=\bigl((3k-2k)/3\bigr)^{1/3}=(k/3)^{1/3}.
\]
Since $\LT(w)\geq \LCS(w_i,w_j)$, the proof is complete.
\end{proof}

There are two natural approaches to improve upon this proof. First, one might hope that
using more than three subsequences would increase the $k^{1/3}$ bound in the Beame and Huynh-Ngoc 
result. Second, instead of replacing copies of a same symbol by different symbols
in \cref{lem:multiperm_simple}, we might hope to take advantage of the repetitions.
We examine these approaches in order.

\paragraph{Common subsequences in large sets of permutations.}
For a set $\W$ of words let 
$\LCS_T(\W)\eqdef \max \LCS(W')$
where the maximum is taken over all sets $\W'\subset\W$ of size $T$. For a family
of words $\mathcal{F}$ let
\[
  \LCS_T(t,\mathcal{F})\eqdef \min_{\W\in \binom{\mathcal{F}}{t}} \LCS_T(\W).
\]
With this notation, the Beame--Huynh-Ngoc theorem asserts that
$\LCS_2(3,\P_k)\geq k^{1/3}$. For a constant number of permutations, a matching
upper bound was proved by Beame, Blais and Huynh-Ngoc \cite[Theorem~2]{bbh}:
\[
  \LCS_2(t,\P_k)\leq 32 (tk)^{1/3}\qquad\text{for all }t\leq k^{1/2}.
\]
We tweak the construction from \cite{bbh} to show optimality
for $t\leq k^{1/3}$:
\begin{theorem}\label{thm:tweak}(Proof is in \cref{sec:constr})
For every $t$ and $k$ satisfying $3\leq t\leq k^{1/2}$, we have
\begin{align}
\LCS_2(t,\P_k)&\leq 4k^{1/3}+O(k^{7/40})&&\text{if }t\leq k^{1/3},\label{eq:tweakone}\\
\LCS_2(t,\P_k)&\leq 4t+O(t^{21/40})&&\text{if }k^{1/3}\leq t\leq k^{1/2}.\label{eq:tweaktwo}
\end{align}
\end{theorem}
We do not know if the constant of $4$ in \cref{thm:tweak} is sharp. However, the next theorem shows that, for many permutations, the constant
is not $1$, as the result of Beame--Huynh-Ngoc might have suggested.
\begin{theorem}\label{thm:lowerbound_perm}(Proof is in \cref{sec:lowerbounds})
For every $t$ we have $\LCS_2(t,\P_k)\geq (\tfrac{7}{4}-\tfrac{8}{t})^{1/6}k^{1/3}-2$.
\end{theorem}

The connections between $\LCS_2$ and $\LT$ go both ways: The next
result, \cref{thm:twinrecurs} translates upper bounds on $\LCS_2$ into upper bounds on $\LT$.
In the opposite direction, we posit a conjecture asserting a lower bound on $\LCS_2$ ``on average'',
and show in \cref{thm:conjectured_ltlowerbounds} that the conjecture implies non-trivial lower bounds on $\LT$.
\begin{theorem}\label{thm:twinrecurs}(Proof is in \cref{proof:twinrecurs})
If there exists a constant $C$ such that for all $k$ we have $\LCS_2(2k,\P_k)\leq Ck^{1/3}$, then for all $k$ we have $\LT(k,n)\leq 6Ck^{-2/3}n+o(n)$.
\end{theorem}
The following conjecture is a generalization of a well-known fact that the expected length of the longest increasing subsequence
in a permutation sampled uniformly from $\P_k$ has length at least $\sqrt{k}$. 
\begin{conjecture}\label{conj:permlcs}
Consider an arbitrary probability distribution on $\P_k$. Let $\pi_1,\pi_2$ be
two permutations sampled independently from $\P_k$. Then 
$\E\bigl[\LCS(\pi_1,\pi_2)\bigr] \geq \sqrt{k}$.  
\end{conjecture}
It might even be true that $\E\bigl[\LCS(\pi_1,\pi_2)]$ is minimized for the uniform measure on $\P_k$; in that case, the
bound would be asymptotic to $2\sqrt{k}$ by the work of Logan--Shepp \cite{logan_shepp} and Ver\v{s}ik--Kerov \cite{versik_kerov}.
It is straightforward to verify this stronger form of the conjecture for $k\leq 3$. The authors
have no idea how to approach the conjecture (or its negation).

The truth of the conjecture implies a lower bound for $\LCS_2(t,\P_k)$ for large $t$. Indeed,
given a set of $t$ permutations, one can apply the conjecture to the uniform probability
distribution on that set. In that case, the probability that $\pi_1=\pi_2$ is $1/t$,
and so we deduce that $\frac{k}{t}+\LCS_2(t,\P_k)\geq \sqrt{k}$.

As promised, the conjecture implies an improvement on \cref{thm:klarge}:
\begin{theorem}\label{thm:conjectured_ltlowerbounds}(Proof is in \cref{sec:conjlower})
\Cref{conj:permlcs} implies that $\LT(k,n)/n\geq \tfrac{1}{100}k^{-12/19}\log^{-8/19} k-\nobreak o(1)$ for any fixed $k$ and large $n$.
\end{theorem}

\paragraph{Common subsequences in small sets of multipermutations.}
As expected, we can improve the bound in \cref{lem:multiperm_simple} by using
letter repetitions:
\begin{theorem}\label{thm:multiperm_good}(Proof is in \cref{sec:lowerbounds})
Let $\pi_1,\pi_2,\pi_3\in \P_{\vec{s}}$ be three $\vec{s}$-multipermutations of
length $\norm{\vec{s}}=n$. Then there is a pair $i<j$ such that
$\LCS(w_i,w_j)\geq \Bigl(\frac{1}{6}\sum_{l\in [k]} s_l^2 \Bigr)^{1/3}$.
In other words, 
$
  \LCS_2(3,\P_{\vec{s}})\geq \Bigl(\tfrac{1}{6}\sum_{l\in [k]} s_l^2 \Bigr)^{1/3}.
$
\end{theorem}
\begin{theorem}\label{thm:multiperm_sharp}(Proof is in \cref{sec:constr})
The previous bound is sharp: If $\vec{s}=(s,\dotsc,s)\in \Z_+^k$ and $s\leq \tfrac{1}{5}k$,
then $\LCS_2(4,\P_{\vec{s}})\leq (2s^2k)^{1/3}+\tfrac{5}{3}s+s^{4/3}k^{-1/3}$.
\end{theorem}

Surprisingly, this improvement does not hold in the
setting of the classical Erd\H{o}s--Szekeres theorem! 

To explain the meaning of the preceding exclamation we must first recall 
the statement of the Erd\H{o}s--Szekeres theorem.
The \emph{reverse} of a word $w$, denoted $\rev w$, is the word obtained by writing $w$ backward.
For example, $\rev \wordf{abracadabra} = \wordf{arbadacarba}$. Let 
\[
 \LCSr(w,w')\eqdef \max\bigl(\LCS(w,w'),\LCS(w,\rev w')\bigr).
\]
The Erd\H{o}s--Szekeres theorem \cite[p.~467]{erdos_szekeres} asserts that if $\pi=\wordf{12\dots k}$,
then $\LCSr(\pi,\pi')\geq \sqrt{k}$ for every $\pi'\in \P_k$. Since $\LCSr(\pi,\pi')$ is unchanged
by relabelling the alphabet, the inequality $\LCSr(\pi,\pi')\geq \sqrt{k}$ holds
for every $\pi,\pi'\in \P_k$. Hence, by the same reasoning
as in the proof of \cref{lem:multiperm_simple} it follows that $\LCSr(\pi,\pi')\geq \sqrt{n}$
for every two $\vec{s}$-multipermutations $\pi,\pi'$ of length $n$. In view of \cref{thm:multiperm_good}, 
it is quite a surprise then that this bound is sharp!
\begin{theorem}\label{thm:es_sharp}(Proof is in \cref{sec:constr})
If $\vec{s}=(s,\dotsc,s)\in \Z_+^k$ and $n=\norm{\vec{s}}$,  then 
there exist two $\vec{s}$-multipermutations
$\pi,\pi'$ 
such that $\LCSr(\pi,\pi')\leq \sqrt{n}+s$. 
\end{theorem}

\paragraph{Common subsequences in large sets of multipermutations.} 
In an attempt to improve \cref{thm:klarge} it is natural to combine the two approaches, and
consider many multipermutations. Alas, we have been unable to extend \cref{thm:multiperm_sharp} 
or to prove better lower bounds on $\LCS_2(t,\P_{\vec{s}})$.

\section{Regularity lemma and proof \texorpdfstring{\cref{thm:kthree}}{theorem 2}}\label{sec:regthree}
The key ingredient in the proof of \cref{thm:app_binary} was a regularity lemma for words.
We state a version of the lemma that we need.

For a word $w\in [k]^n$ and another word $u$ of length $l$ that is smaller than $n$,
we define \emph{frequency of $u$ in $w$} to be the probability that a randomly
chosen $l$-letter-long subword of $w$ is a copy of $u$. We denote the frequency by
$f_w(u)$. A word $w\in [k]^n$ is \emph{$(\veps,L)$-regular} if whenever
$w'$ is a subword of $w$ of length at least $\veps n$, then
\[
  \bigl\lvert f_w(u)-f_{w'}(u)\bigr\rvert<\veps\qquad\text{for every }u\text{ of length at most }L.
\]
\begin{lemma}[Regularity lemma for words]\label{lem:reg}
For every $\veps>0$ and every $L$ there is a number $M=M(\veps,k,L)$ such that the following holds: Every 
sufficiently long word $w\in [k]^n$ can be partitioned into at most $M$ subwords such that the total length of the subwords that are not
$(\veps,L)$-regular is at most~$\veps n$.
\end{lemma}
This lemma is slightly different than what appears in~\cite[Theorem~6]{app_twins}. The difference is that
in~\cite{app_twins} the result was asserted only for $L=1$. However,   
the general case can be deduced from this special case:
\begin{proof}[Proof that the special case of $L=1$ implies the case of general $L$]
Assume $k,\veps$ and $L$ are fixed. Let $\mathcal{W}$ be the set of all words in $[k]$ of length at most
$L$. We claim that we may take $M(\veps,k,L)=M(\veps/2,2^{\abs{\mathcal{W}}},1)$. Indeed, given a word
$w\in [k]^n$ we define a word $W$ in the alphabet $2^{\mathcal{W}}$ via
\[
  W[i]\eqdef \{ u\in \mathcal{W} : \text{the subword starting from }w[i]\text{ of length }\operatorname{len}(u)\text{ is equal to }u\}.  
\]
Since a subword of $W$ that is
$(\veps/2,1)$-regular corresponds to a $(\veps,L)$-regular subword of $w$, the result follows
from the $L=1$ case of the regularity lemma applied to $W$. (The reason for $\veps/2$ deteriorating into
$\veps$ is the edge effect --- the words that start too close to a subword boundary are miscounted.)
\end{proof}

For brevity, we call a word $w$ simply \emph{$\veps$-regular} if it is $(\veps,1/\veps)$-regular.

In context of proving the lower bounds on the twin word problem, the regularity lemma allows us to
assume that the word under consideration is $\veps$-regular, for any fixed $\veps$. Indeed, suppose the
bound of $\LT(u)\geq \alpha \operatorname{len}(u)$ is valid for all $\veps$-regular words $u$ of length exceeding $n_0$, and $w=w_1\dotsb w_m$
is the partition into $m\leq M$ parts described in the lemma. Then the bound $\LT(u)\geq \alpha\bigl(\operatorname{len}(u)-n_0\bigr)$ is
valid for all $\veps$-regular words, and so $\LT(w)\geq \alpha (n-mn_0)-\veps n$.

We shall prove \cref{thm:kthree} for $k=3$. For a general $k$ \cref{thm:kthree} follows from the $k=3$
case by stripping all but the three most popular letters from a word. So, we assume that $w\in [3]^n$
is $\veps^2$-regular for some fixed, but arbitrarily small $\veps>0$. Set
\[
  \beta\eqdef \min\bigl(f(\wordf{1}),f(\wordf{2}),f(\wordf{3})\bigr).
\]
We may also assume that
$\beta >0.02$, for otherwise
$\LT(w)\geq 0.49n-o(n)$ follows from \cref{thm:app_binary} applied to the two most-frequent letters in $w$.

For brevity, we write $f(u)$ in place
of $f_w(u)$ throughout the remainder of this section, and also $f(u_1 + u_2 + \dotsb)$ in place of
$f(u_1)+f(u_2)+\dotsb$.

\begin{lemma}\label{lem:runbound}
We have $\LT(w)/n\geq \tfrac{1}{2}f(\wordf{11}+\wordf{22}+\wordf{33})-O(1/n)$.
\end{lemma}
\begin{proof}
Write $w=w_1w_2\dotsb w_t$ be a partition of $w$ into subwords that consist of a single letter of the alphabet.
For example, $\wordf{2223312222111}$ would be partitioned as $\wordf{222}\ \wordf{33}\ \wordf{1}\ \wordf{2222}\ \wordf{111}$. The number
of subwords is equal to
\[
  1+(n-1)f(\wordf{12}+\wordf{13}+\wordf{21}+\wordf{23}+\wordf{31}+\wordf{32})=n\bigl(1-f(\wordf{11}+\wordf{22}+\wordf{33})\bigr)+O(1).
\]
Since $\LT(w_i)\geq \operatorname{len}(w_i)/2-1/2$, and $\LT(w)\geq \sum_i \LT(w_i)$, the lemma follows from $\sum_i \operatorname{len}(w_i)=n$.
\end{proof}
\begin{lemma}\label{lem:cleverpairing}
Define
$\alpha_1=f(\wordf{21}+\wordf{31})$,
$\alpha_2=f(\wordf{12}+\wordf{32})$, 
$\alpha_3=f(\wordf{13}+\wordf{23})$. 
For each $l\in [3]$ we then have 
\[
  \LT(w)/n\geq \tfrac{1}{2}\Bigl(1-f(l)+\frac{\alpha_l^2}{1-f(l)}\Bigr)-O(\veps).
\]
\end{lemma}
\begin{proof}
In view of the symmetry it suffices to prove the case $l=\wordf{3}$.
Let $t\eqdef \lfloor 1/\veps\rfloor$. Partition $w$ into $t$ subwords of 
length at least $\veps n$ each; say, $w=w_1w_2\dotsb w_t$. 
Note that since $w$ is  $\veps^2$-regular, each of $w_i$ is $\veps$-regular.

We will find twins of the following form:\smallskip
\begin{center}
\begin{tabular}{rccccccc}
           & $w_1$                       & $w_2$                       & $w_3$ & $\dotsb$&$w_{t-1}$&$w_t$\\\cline{1-7}
First twin & \wordf{1}'s and \wordf{3}'s & \wordf{2}'s and \wordf{3}'s & \wordf{1}'s and \wordf{3}'s & $\dotsb$ & \wordf{2}'s and \wordf{3}'s\\
Second twin &                            & \wordf{1}'s and \wordf{3}'s & \wordf{2}'s and \wordf{3}'s & $\dotsb$ & \wordf{1}'s and \wordf{3}'s & \wordf{2}'s and \wordf{3}'s
\end{tabular}
\end{center}\medskip

So, the first twin will contain \wordf{1}'s and \wordf{3}'s, but no \wordf{2}'s from $w_1$, etc. To assure that
no letter is in both twins, we adopt the following rule: Only a \wordf{3} that immediately
follows a \wordf{1} (for $w_1,w_3,\dotsc$) or a~\wordf{2} (for $w_2,w_4,\dotsc$) in $w$
can appear in the first twin. Similarly, only a~\wordf{3} that immediately follows
a~\wordf{2} (for $w_3,w_5,\dotsc$) or a~\wordf{1} (for $w_2,w_4,\dotsc$) in $w$
can appear in the second twin. 

Next, we show how to find a long common subsequence between $w_1$ and $w_2$ 
consisting only of \wordf{1}'s and \wordf{3}'s that satisfies the restriction
on \wordf{3}'s specified above. We will do so by first finding subsequences
$u_1$ and $u_2$ of $w_1$ and $w_2$ respectively that consist only of \wordf{1}'s,
and then adding \wordf{3}'s where possible.

Let $r_1$ and $r_2$ be the number of \wordf{1}'s in $w_1$ and $w_2$, respectively.
Set $r\eqdef \min(r_1,r_2)-\veps^2 n$. Note that $r=\bigl(f(\wordf{1})-O(\veps)\bigr)\veps n$
by the regularity of $w$. Let $u_1$ be the subsequence of $w_1$ that consists of the 
first $r$ occurrences of \wordf{1} in $w_1$.  Suppose that the \wordf{1}'s in $w_2$
are at positions $i_1,i_2,\dotsc$. Pick an integer $m$ uniformly at random 
from $0$ to $\veps^2 n$, and let $u_2$ be the subsequence 
$w_2[i_m], w_2[i_{m+1}],\dotsc,w_2[i_{m+r-1}]$. 

As words, $u_1$ and $u_2$ are both words of length $r$ that contain only \wordf{1}'s.
As subsequences of $w$, they are more interesting. Of $r$ letters that $u_1$
contains, $\bigl(f(\wordf{13})-O(\veps)\bigr)\veps n$ are followed by a \wordf{3} in $w_1$.
Say, $u_1[i]$ is followed by a \wordf{3} in $w_1$. Consider the
letter $u_2[i]$. Due to the choice of $m$, the $u_2[i]$ is chosen
uniformly from all \wordf{1}'s in an interval of length at least $\veps^2 n$.
So, crucially, since $w_2$ is $\veps$-regular, the probability that
$u_2[i]$ is followed by a \wordf{3} is $\frac{f(\wordf{13})-O(\veps)}{f(\wordf{1})-O(\veps)}=f(\wordf{13})/f(\wordf{1})-O(\veps)$. By linearity of expectation, this implies that there is an $m$
such that for at least
\[
  \frac{f(\wordf{13})-O(\veps)}{f(\wordf{1})}\bigl(f(\wordf{13})-O(\veps)\bigr)\veps n=
  \frac{f(\wordf{13})^2-O(\veps)}{f(\wordf{1})}\veps n
\]
values of $i$ both $u_1[i]$ and $u_2[i]$ are followed by a \wordf{3}. Hence, we can
extend $u_1$ and $u_2$ to subsequences $u_1'$ and $u_2'$ of $w_1$ and $w_2$, respectively,
of length at least
\[
  \left(f(\wordf{1})+\frac{f(\wordf{13})^2-O(\veps)}{f(\wordf{1})}\right)\veps n.
\]

Similar matches can be found between $w_2$ and $w_3$, between $w_3$ and $w_4$, etc. 
Concatenation of these matches yields a pair of twins that are large:
\begin{align*}
  \LT(w)&\geq 
  \lfloor t/2\rfloor\left(f(\wordf{1})+\frac{f(\wordf{13})^2-O(\veps)}{f(\wordf{1})}\right)\veps n+  
  \lfloor (t-1)/2\rfloor\left(f(\wordf{2})+\frac{f(\wordf{23})^2-O(\veps)}{f(\wordf{2})}\right)\veps n\\
  &=\tfrac{1}{2}\left(f(\wordf{1})+f(\wordf{2})+\frac{f(\wordf{13})^2}{f(\wordf{1})}+\frac{f(\wordf{23})^2}{f(\wordf{2})}-O(\veps)\right)n,
  \\\intertext{which by the Cauchy--Schwarz inequality applied to the vectors $\bigl(\frac{f(\wordf{13})}{\sqrt{f(\wordf{1})}},\frac{f(\wordf{23})}{\sqrt{f(\wordf{2})}}\bigr)$ and $\bigl(\sqrt{f(\wordf{1})},\sqrt{f(\wordf{2})}\bigr)$ implies that}
  2\LT(w)/n&\geq f(\wordf{1})+f(\wordf{2})+\frac{\alpha_3^2}{f(\wordf{1})+f(\wordf{2})}-O\bigl(\veps\bigr).
\end{align*}
Since $f(\wordf{1})+f(\wordf{2})=1-f(\wordf{3})$, the proof of the lemma is complete.
\end{proof}

The preceding two lemmas are enough to deduce \cref{thm:kthree}. 
Indeed, applying \cref{lem:cleverpairing} for $l=1,2,3$
and adding the resulting bounds we obtain
\begin{align*}
  6\LT(w)/n&\geq 3-f(\wordf{1}+\wordf{2}+\wordf{3})+\frac{\alpha_1^2}{1-f(\wordf{1})}+\frac{\alpha_2^2}{1-f(\wordf{2})}+\frac{\alpha_3^2}{1-f(\wordf{3})}-O(\veps)\\
  &= 2+\frac{\alpha_1^2}{1-f(\wordf{1})}+\frac{\alpha_2^2}{1-f(\wordf{2})}+\frac{\alpha_3^2}{1-f(\wordf{3})}-O(\veps)\\
  &\geq 2+\tfrac{1}{2}(\alpha_1+\alpha_2+\alpha_3)^2-O(\veps),
\end{align*}
where the last line follows from applying Cauchy--Schwarz inequality to vectors $\bigl(\frac{\alpha_1}{\sqrt{1-f(\wordf{1})}},\frac{\alpha_2}{\sqrt{1-f(\wordf{2})}},\frac{\alpha_3}{\sqrt{1-f(\wordf{3})}}\bigr)$ and $\bigl(\sqrt{1-f(\wordf{1})},\sqrt{1-f(\wordf{2})},\sqrt{1-f(\wordf{3})}\bigr)$, and then using $f(\wordf{1}+\wordf{2}+\wordf{3})=1$ to simplify the resulting
expression.
Since $\alpha_1+\alpha_2+\alpha_3=1-f(\wordf{11}+\wordf{22}+\wordf{33})$, 
in view of the bound from \cref{lem:runbound} we conclude that
\[
  \LT(w)/n\geq \min_x \max (\tfrac{1}{3}+\tfrac{1}{12}x^2,\tfrac{1}{2}-\tfrac{1}{2}x)-O(\veps)=\frac{4-\sqrt{11}}{2}-O(\veps).
\]
Since $\frac{4-\sqrt{11}}{2}>\tfrac{1}{3}\cdot 1.02$ and $\veps$ is arbitrary, \cref{thm:kthree} follows.

\section{Proof of \texorpdfstring{\protect\cref{thm:random_bound}}{theorem 4}}\label{sec:random_bound}
We will show that with high probability a random word of length $n$ satisfies the conclusion
of \cref{thm:random_bound}.  Recall that $w[i]$ denotes the $i$'th
letter of the word $w$. 

Twins $w_1$ and $w_2$ in $w\in [k]^n$ are said to be \emph{monotone} if 
$w_1[i]$ precedes $w_2[i]$ in $w$ for all $i$.
\begin{lemma}
If $\LT(w)\geq m$, then $w$ contains monotone twins of length $m$.
\end{lemma}
\begin{proof}
The condition implies that $w$ contains twins of length $m$. 
However, if $w[p_1]\dotsb w[p_m]$ and $w[p_1']\dotsb w[p_m']$ are twins, then so
are $w[\bar{p}_1]\dotsb w[\bar{p}_m]$ and $w[\bar{p}_1']\dotsb w[\bar{p}_m']$,
where $\bar{p}_i=\min(p_i,p_i')$ and $\bar{p}_i'=\max(p_i,p_i')$.
\end{proof}

To each pair $(w_1,w_2)$ of monotone twins in $w$ we associate a word
$R(w_1,w_2,w)\in \{\wordf{0},\wordf{1},\wordf{2}\}^n$ by the rule
\[
  R(w_1,w_2,w)[i]=
   \begin{cases}
     \wordf{0}&\text{if }w[i]\text{ is neither in }w_1\text{ nor in }w_2,\\
     \wordf{1}&\text{if }w[i]\text{ is in }w_1,\\
     \wordf{2}&\text{if }w[i]\text{ is in }w_2.\\
   \end{cases}
\]
The word $R(w_1,w_2,w)$ records the \emph{roles} of letters of $w$ in the pair $(w_1,w_2)$.
For example, consider the monotone pair 
$\wordf{1}\overline{\wordf{0}}\underline{\wordf{0}}\overline{\wordf{11}}\underline{\wordf{1}}\wordf{0}\underline{\wordf{1}}\overline{\wordf{1}}\underline{\wordf{1}}\wordf{01}$,
where the overlines indicate the letters in $w_1$ and the underlines indicate the letters in $w_2$. For this pair, $R=\wordf{012112021200}$.

A monotone pair $(w_1,w_2)$ in $w$ is said to be \emph{regular} if the following two conditions hold:
\begin{enumerate}
\item There exist no two numbers $i<j$ satisfying the following: $R[i]=\wordf{2}$ and $R[j]=\wordf{1}$, and 
$R[k]=\wordf{0}$ for all $i<k<j$, and $w[i]=w[j]$.
\item There exist no two numbers $i<j$ satisfying the following: $R[i]\in\{\wordf{1},\wordf{2}\}$ and $R[j]= \wordf{0}$, and 
$R[k]=\wordf{0}$ for all $i<k<j$, and $w[i]=w[j]$.
\end{enumerate}
We can express these conditions in the overline/underline notation: The condition (a) forbids the pattern $\underline{\wordf{x}}\wordf{???}\overline{\wordf{x}}$,
whereas the condition (b) forbids the patterns $\overline{\wordf{x}}\wordf{???x}$
and $\underline{\wordf{x}}\wordf{???x}$, where the question marks denote letters that are not
in the twins.

\begin{lemma}\label{lem:regtwin}
If $\LT(w)\geq m$, then $w$ contains regular twins of length $m$. 
\end{lemma}
\begin{proof}
Pick monotone twins $(w_1,w_2)$ in $w$ such that $R(w_1,w_2,w)$ is lexicographically minimal. Then $(w_1,w_2)$
is regular. Indeed, if $(w_1,w_2)$ were not regular, then swapping the roles of $w[i]$ and $w[j]$
would lead to monotone twins with a lexicographically smaller value of $R$.
\end{proof}

Let $\mathcal{R}_m^n$ consist of all words $R \in \{\wordf{0},\wordf{1},\wordf{2}\}^n$ in which letter $\wordf{1}$ and letter $\wordf{2}$
occur $m$ times each. For $R\in \{\wordf{0},\wordf{1},\wordf{2}\}^n$ let $p(R)$ be the number of occurrences of the
pattern $\wordf{20}^{*}\wordf{1}$, i.e., a $\wordf{2}$ followed by zero or more $\wordf{0}$'s, and then
followed by a $\wordf{1}$. Also, let $z(R)$ be the length of the longest prefix of $R$ that contains only $\wordf{0}$'s.
Let $\mathcal{R}_{m,p,z}^n\eqdef \{R \in \mathcal{R}_m^n : p(R)=p,\ z(R)=z\}$.
For example, $\wordf{012112021200}\in \mathcal{R}_{4,2,1}^{12}$.
 
Let $\mathcal{M}\subset \{\wordf{0},\wordf{1},\wordf{2}\}^n$ be the set of all words in which every prefix contains at least as many \wordf{1}'s
as \wordf{2}'s. Note $R(w_1,w_2,w)\in \mathcal{M}$ for every pair of monotone twins. 
For $R\in \mathcal{M}$ let $B_R$ be the event that a word $w$ chosen uniformly
at random from $[k]^n$ contains a regular pair $(w_1,w_2)$ satisfying
$R(w_1,w_2,w)=R$.
\begin{lemma}\label{lem:regtwinprob}
Suppose $R\in \mathcal{R}_{m,p,z}^n\cap \mathcal{M}$, then $\Pr[B_R]=(1/k)^m (1-1/k)^{p+n-2m-z}$.
\end{lemma}
\begin{proof}
Imagine that each letter $w[i]$ of $w$ is in its own box, which is labeled $R[i]$. Call a box labeled with \wordf{0} a \emph{\wordf{0}-box}, etc. 
Imagine also that boxes are connected by red and blue wires according to the following rules:
\begin{itemize}
\item The $i$'th box labeled \wordf{1} is connected by a blue wire to the $i$'th box labeled \wordf{2}.
\item If the first non-zero-labeled box that precedes a \wordf{1}-box is a \wordf{2}-box,
then  the \wordf{1}-box and the \wordf{2}-box are connected by a red wire.
\item Each \wordf{0}-box is connected with the preceding non-zero-labeled box by a red wire.
\end{itemize}

The following condition is clear from the definition of a regular pair:\medskip
\begin{center}
\parbox{30em}{%
A word $w$ contains a regular pair $(w_1,w_2)$ satisfying
$R(w_1,w_2,w)=R$ if and only if the blue wires connect the boxes containing the same letters,
and red wires connect the boxes containing different letters.}
\end{center}\medskip

The boxes start closed, and we open them one by one. Each time
we open a box we look at the letter that it contains, and at the letters of all 
previously-opened boxes that this box is connected to. If the condition above is violated,
we abort. 

The order for opening the boxes is not arbitrary: we first open
all non-zero boxes from left to right, and only then open the \wordf{0}-boxes.
This order ensures that each time we open a box, there is at most one wire
that connects it to a previously-opened box. Thus, the probability that
we abort is $1-1/k$ if that the wire is blue, and $1/k$ if the wire is red.

Since the number of blue wires is $m$, and the number of red wires
is $p+(n-2m-z)$, the lemma follows.
\end{proof}

\begin{lemma}
The size of $\mathcal{R}_{m,p,z}^n$ is 
$\abs{\mathcal{R}_{m,p,z}^n}=\binom{n-z}{2m}\binom{m}{p}^2$.
\end{lemma}
\begin{proof}
To each word $R\in \mathcal{R}_{m,p,z}^n$ associate the
word $R'\in \mathcal{R}_{m,p,0}^{2m}$ obtained by removing all \wordf{0}'s from $R$.
Since $R$ starts with a prefix of $z$ \wordf{0}'s, one can recover $R$
from $R'$ by specifying the positions of the remaining \wordf{0}'s.
So, $\abs{\mathcal{R}_{m,p,z}^n}=\binom{n-z}{2m}\abs{\mathcal{R}_{m,p,0}^{2m}}$,
and to complete the proof it suffices to compute $\abs{\mathcal{R}_{m,p,0}^{2m}}$.

Every $R\in \mathcal{R}_{m,p,0}^{2m}$ is necessarily of the form
\[
  \fbox{\phantom{123}}\,\wordf{21}\,\fbox{\phantom{123}}\,\wordf{21}\,\fbox{\phantom{123}}\,\wordf{21}\,\dotsb\,\wordf{21}\,\fbox{\phantom{123}}\,\wordf{21}\,\fbox{\phantom{123}}\,,
\]
where there are $p+1$ bins separated by the $p$ occurrences of $\wordf{21}$,
and each bin contains a subword of the form $\wordf{1}^*\wordf{2}^*$ (the subword might be empty).
The bins contain a total of $m-p$ \wordf{1}'s and the same
number of \wordf{2}'s. 
There are $\binom{(m-p)+(p+1)-1}{(p+1)-1}=\binom{m}{p}$
ways of placing $m-p$ identical objects into $p+1$ labeled bins. In particular,
there are $\binom{m}{p}$ ways to place $m-p$ \wordf{0}'s into $p+1$ bins,
and the same number of ways to place $m-p$ \wordf{1}'s into the bins.
Since the choices for placement of \wordf{0}'s and \wordf{1}'s are independent,
we conclude that 
$\abs{\mathcal{R}_{m,p,0}^{2m}}=\binom{m}{p}^2$ as promised.
\end{proof}

Let $m=\alpha n$ for the constant $\alpha$ from \cref{thm:random_bound}.
The union bound and the three preceding lemmas imply that
\begin{align*}
  \Pr[\LT(w)\geq m]&\leq \sum_{p,z} \sum_{R\in\mathcal{R}_{m,p,z}^n\cap \mathcal{M}} \Pr[B_R]\leq
  \sum_{p,z\geq 0} \binom{n-z}{2m}\binom{m}{p}^2 (1/k)^m (1-1/k)^{p+n-2m-z}\\
  &=(1/k)^m(1-1/k)^{n-2m} \sum_{p} \binom{m}{p}^2 (1-1/k)^p \sum_{z\geq 0} \binom{n-z}{2m} (1-1/k)^{-z} \\
  &\leq (1/k)^m(1-1/k)^{n-2m} n^2 \max_p \binom{m}{p}^2 (1-1/k)^p \max_{z\geq 0} \binom{n-z}{2m} (1-1/k)^{-z}.
\end{align*}
Let $f(p)=\binom{m}{p}^2 (1-1/k)^p$ and $g(z)=\binom{n-z}{2m} (1-1/k)^{-z}$. 
Since
$g(z+1)/g(z)=\frac{n-z-2m}{n-z}\cdot \frac{k}{k-1}$ and $m\geq n/2k$ 
the maximum of $g(z)$ is attained at $z=0$. Similarly,
$f(p+1)/f(p)=\left(\frac{m-p}{p+1}\right)^2 (1-1/k)$ implies
that the maximum of $f(p)$ is attained when $p=\frac{m}{1+(1-1/k)^{-1/2}}+O(1)$.
We plug these into the displayed formula above, use
the asymptotic formula $\frac{\log \binom{n}{\beta n}}{n}=\beta \log \frac{1}{\beta}+(1-\beta)\log \frac{1}{1-\beta}+o(1)$, and simplify to obtain:
\begin{equation*}
  \frac{\log \Pr[\LT(w)\geq \alpha n]}{n}\leq (1-2\alpha)\log\frac{1}{1-2\alpha}-\alpha\log(\alpha^2k)-2\alpha\log \left(\frac{2}{1+\sqrt{1-1/k}}\right)+(1-2\alpha)\log(1-1/k) +o(1).
\end{equation*}
Whenever the expression on the right is negative, we have $\Pr[\LT(w)\geq \alpha n]<1$
for large enough $n$, and so $\LT(n)<\alpha n$ for those $n$.

We note that the first two terms in the inequality above are the same 
as in the bound obtained in \cite{app_twins}. The last two terms are new to our analysis.

The bounds on $\LT(4,n)$, $\LT(5,n)$ and the asymptotic bound on
$\LT(k,n)$ for large $k$ in \cref{thm:random_bound} were obtained using the   
\textsc{Mathematica} software package. The code that we used is available at
\url{http://www.borisbukh.org/code/twins_lcs13.html}.

\section{Proof of \texorpdfstring{\cref{thm:twinrecurs}}{theorem 8}}\label{proof:twinrecurs}
\begin{lemma}
Suppose $k,K,A,B$ are natural numbers that satisfy $\LT(K,n)\leq A$ and $\LCS_2(K,\P_k)\leq B$. Then
$\LT(k,kn)\leq (2n-1)B+kA$.
\end{lemma}
\begin{proof}
Let $w\in [K]^n$ be a word such that $\LT(w)\leq A$. Let
$\pi_1,\dotsc,\pi_K\in \P_k$ be a family of $K$ permutations
such that $\LCS(\pi_i,\pi_j)\leq B$ for all distinct $i,j$.

Replace each letter $l$ in $w$ by the permutation $\pi_l$
to obtain a word $w'\in [k]^{kn}$. We claim that $\LT(w')\leq (2n-1)B+kA$.
Indeed, suppose $w_1'$ and $w_2'$ are twins in $w'$. 
Let $\mathcal{T}$ be the set of all pairs $(i',j')$ such that
$w'[i']$ is a letter in $w_1'$ that is matched to $w'[j']$,
which is in $w_2'$. Let $\mathcal{H}\eqdef \bigl\{ (\lfloor i'/k\rfloor,\lfloor j'/k\rfloor) : (i',j')\in \mathcal{T} \bigr\}$.
Let $\mathcal{M}\eqdef \{ (i,j) \in \mathcal{H} : w[i]=w[j]\}$. We may also suppose
that the twins $w_1',w_2'$ were chosen among all twins of their length
so that $\sum_{i,j\in \mathcal{H}} (i+j)$ is minimized.

A simple picture explains the meaning of $\mathcal{T}$ and $\mathcal{H}$
just defined. Imagine a $kn$-by-$kn$ square
that is subdivided into $n^2$ $k$-by-$k$ squares, and imagine that the $kn$ letters of $w'$
are laid along each of the axes. The intervals of length $k$ on the axes correspond to the original letters of $w$. 
For each pair of matching letters of $w_1'$ and $w_2'$ draw a corresponding point --- these
will be points of $\mathcal{T}$. The set of $k$-by-$k$ squares that are
hit by $\mathcal{T}$ is the set $\mathcal{H}$. The squares in $\mathcal{H}\setminus \mathcal{M}$
cannot contain more than $B$ points each since they correspond to different
permutations. The points of $\mathcal{T}$ form a graph of a monotone
function, and so $\abs{\mathcal{H}}\leq 2n-1$. Hence, at most
$(2n-1)B$ points of $\mathcal{T}$ fall into $\mathcal{H}\setminus \mathcal{M}$.

It remains to bound $\abs{\mathcal{M}}$. We claim that no two squares in $\mathcal{M}$
share an $x$-coordinate, and no two squares share a $y$-coordinate. 
Indeed, suppose $(i,j_1),(i,j_2)\in \mathcal{M}$ with $j_1<j_2$.
Then by moving all the points from the square indexed by $(i,j_2)$ to
the one indexed by $(i,j_1)$ we would obtain a pair of twins
of same length, but with the smaller value of $\sum_{i,j\in \mathcal{H}} (i+j)$. 
Since this contradicts the minimality assumption, the claim follows.
Put $w_1=\bigl\{ w[i] : (i,j)\in\mathcal{M}\bigr\}$ and $w_2=\bigl\{ w[j] : (i,j)\in\mathcal{M}\bigr\}$.
The preceding claim implies that $w_1$ and $w_2$ are twins in $w$, and so $\abs{\mathcal{M}}\leq A$.
Hence, the number of points of $\mathcal{T}$ in $\mathcal{M}$ is at most $kA$.
Together with the estimate on points $\mathcal{T}$ in $\mathcal{H}\setminus {M}$,
this implies that $\abs{\mathcal{T}}\leq (2n-1)B+kA$. Since $\LT(w')=\abs{\mathcal{T}}$,
the proof is complete.
\end{proof}

Suppose the condition of \cref{thm:twinrecurs} holds, and $C$ is a constant
as in the theorem. Let $\alpha_k$ be the least real number such that
$\LT(k,n)\leq \alpha_k n+o(n)$. The preceding lemma implies that
$\LT(k,kn)\leq 2C n k^{1/3}+k\LT(2k,n)$, and so 
$\alpha_k\leq 2C k^{-2/3}+\alpha_{2k}$. Thus,
\[
  \alpha_k\leq 2C k^{-2/3}+2C (2k)^{-2/3}+\dotsb+2C (2^{t-1} k)^{-2/3}+\alpha_{2^t k}.
\]
As $t\to \infty$, the last term tends to $0$ by inequality \eqref{eq:goestozero}, and the sum of the 
remaining terms tends to $\frac{2C}{1-2^{-2/3}}k^{-2/3}\leq 6Ck^{-2/3}$.

\section{Proof of \texorpdfstring{\cref{thm:conjectured_ltlowerbounds}}{theorem 10}}\label{sec:conjlower}
Throughout this section we employ the following notation. For a word $w$ and a letter $l$ 
we write $l\in w$ if the letter $l$ occurs in $w$.

\begin{lemma}\label{lem:concatlcs}
Suppose $\bar{u}=u_1\dotsb u_r$ and $\bar{u}'=u_1'\dotsb u_r'$ are words that are the 
concatenation of $r$ other words. Then
$\LCS(\bar{u},\bar{u}')\leq \sum_{i,j} \LCS(u_i,u_j')$.
\end{lemma}
\begin{proof} A common sequence of $\bar{u}$ and $\bar{u}'$ can be broken up as a concatenation of common
sequences of $u_i$ and $u_j$ over various $i,j$.\end{proof}

\begin{proof}[Proof of \cref{thm:conjectured_ltlowerbounds}]
In view of the discussion following the regularity lemma in \cref{sec:regthree}, it suffices 
to prove $\LT(w)/n\geq \tfrac{1}{100}k^{-12/19}\log^{-8/9}k-o(1)$ for all  $\veps$-regular words $w\in [k]^n$.
So, we assume that $w\in [k]^n$ is given, and that it is $\veps$-regular. 

Let $m=n/12k$.
Pick an integer $r$ from the interval $[0,n/4]$ uniformly at random, and pick
another integer $r'$ from the interval $[n/2,3n/4]$ uniformly at random independently from $r$. 
Starting from the position $r$ partition $w$ into $m$ intervals of length $3k$ each.
Note that the total length of these intervals is $n/4$, and so these intervals are
completely contained in the first half of $w$. Let $w_1,\dotsc,w_m\in [k]^{3k}$ be the subwords
in these intervals. Similarly, starting from the position $r'$ define $m$ subwords
$w_1',\dotsc,w_m'\in [k]^{3k}$ that are completely contained in the second half of $w$.

Since $w$ is $\veps$-regular for every word $u\in [k]^{3k}$ and every
$i\in [m]$ we have $\Pr[w_i=u]=\mu(\{u\})+O(\veps)$, and $\Pr[w_i'=u]=\mu(\{u\})+O(\veps)$.
Here, $\mu$ is a probability measure on $[k]^{3k}$. Furthermore, $w_i$ and $w_i'$ are independent.
Note that 
\begin{equation}
\label{eq:LTlower}
\LT(w)\geq \sum_{i=1}^m \LT(w_i)\qquad\text{and}\qquad
\LT(w)\geq \sum_{i=1}^m \LCS(w_i,w_i'). 
\end{equation}
We shall show that at least one of these two bounds is large.

Let $0<\alpha<1/4$ be a parameter to be chosen later.  Put 
\[
  S\eqdef \{ u\in [k]^{3k} : \text{ fewer than }\alpha k\text{ distinct letters occur in }u\}.
\]

We distinguish two cases, depending on $\mu(S)$. 
\begin{enumerate}
\item Suppose $\mu(S)\geq \tfrac{1}{2}$. If $u\in S$ is arbitrary, then treating $u$ as a word over an $\alpha k$-letter alphabet, and using \cref{thm:klarge}
we conclude that
\[
  \LT(u)\geq 3^{-4/3} (\alpha k)^{-2/3}3k-3^{-1/3}(\alpha k)^{1/3}=3^{-1/3}\alpha^{-2/3}(1-\alpha)k^{1/3}\geq \tfrac{1}{2}\alpha^{-2/3}k^{1/3}.
\]
Thus, $\LT(w)\geq \E\bigl[\sum_i \LT(w_i)\bigr]\geq m\cdot \bigl(\mu(S)-O(\veps)\bigr)\tfrac{1}{2}\alpha^{-2/3}k^{1/3}\geq \tfrac{1}{48}\alpha^{-2/3}k^{-2/3}-O(\veps n).$

\item Suppose $\mu(\bar{S})\geq \tfrac{1}{2}$.
Let $u$ and $u'$ be two words sampled independently from $[k]^{3k}$ according to measure $\mu$.
Let $L=\{ l \in [k] : \Pr[l\in u]\geq \alpha/4\}$. Since $\mu(\bar{S})\geq \tfrac{1}{2}$,
for random $l\in[k]$ we have $\Pr_{l,u} [l\in u]\geq \alpha/2$ and so
it follows that $\abs{L}\geq \alpha k/4$. For a word $v$ that contains every
letter of $L$ at least once denote by $\pi_L(v)$ the permutation on the 
alphabet $L$ obtained by taking in $v$ the first occurrence of each letter.
If there is a letter of $L$ that does not occur in $v$, we let $\pi_L(v)=\emptyword$,
where $\emptyword$ is the empty word.

Let $r\eqdef 4\alpha^{-1}\log (4k)$. Let $u_1,\dotsc,u_r$ and $u_1',\dotsc,u_r'$ be
$2r$ words sampled independently from $[k]^{3k}$ according to measure $\mu$. Let $\bar{u}=u_1\dotsb u_r$
and $\bar{u}'=u_1'\dotsb u_r'$. Put $\pi=\pi_L(\bar{u})$ and $\pi'=\pi_L(\bar{u}')$. 
For each fixed $l\in L$, we have 
$\Pr[l\not\in \bar{u}]\leq (1-\alpha/4)^r\leq 1/4k$. Hence, 
$\Pr[\pi=\emptyword]\leq 1/4$. Similarly $\Pr[\pi'=\emptyword]\leq 1/4$. 
By \cref{conj:permlcs} it follows
that 
\begin{align*}
  \E\bigl[\LCS(\pi,\pi')\bigr]&=\E\bigl[\LCS(\pi,\pi') | \pi\neq\emptyword\,\wedge\,\pi'\neq \emptyword \bigr]\Pr[\pi\neq\emptyword\,\wedge\,\pi'\neq \emptyword]\\
  &\geq\tfrac{1}{2}\sqrt{\abs{L}}\geq \tfrac{1}{4}\sqrt{\alpha k}.
\end{align*}
By \cref{lem:concatlcs} and the linearity of expectation we have
\[
  r^2\E\bigl[\LCS(u,u')\bigr]=\sum_{i,j\in [r]} \E\bigl[\LCS(u_i,u_j')\bigr]\geq \E\bigl[\LCS(\bar{u},\bar{u}')\bigr].
\]
Since $\LCS(\bar{u},\bar{u}')\geq \LCS(\pi,\pi')$, we can combine
the two preceding inequalities with \eqref{eq:LTlower} to obtain
\[
\LT(w)\geq \E\Bigl[\sum_{i=1}^m \LCS(w_i,w_i')\Bigr]\geq m\cdot \bigl(\tfrac{1}{r^2}\cdot\tfrac{1}{4}\sqrt{\alpha k}-O(\veps)\bigr)\geq \tfrac{1}{1000}\alpha^{5/2}\frac{k^{-1/2}}{\log^2 k}n -O(\veps n).
\]
\end{enumerate}
Thus no matter which of the two cases holds, we have
\[
  \frac{\LT(w)}{n}\geq \min\bigl(\tfrac{1}{48}\alpha^{-2/3}k^{-2/3},\tfrac{1}{1000}\alpha^{5/2}\tfrac{k^{-1/2}}{\log^2 k}\bigr)-O(\veps).
\]
Since $\veps>0$ is arbitrary, setting $\alpha=3k^{-1/19}\log^{12/19}k$ yields the desired result.
\end{proof}

\section{Constructions}\label{sec:constr}
In this section we prove \cref{thm:tweak,thm:es_sharp,thm:multiperm_sharp}. All of these results are proved by exhibiting explicit constructions
of families of (multi)permutations with required properties. The constructions are very similar to
one another, and so we start by describing their commonality.

In all the constructions the alphabet is $X\times Y\times \dotsb$ 
for some sets $X,Y,\dotsc$. This way each letter $l$ can be thought 
of as a vector $(x,y,\dotsc)$. Given an injective map 
$f\colon X\times Y\times \dotsb \to \Z\times \dotsb \times \Z$
we define $\pi_f$ to be the permutation in which letter $l$ precedes
letter $l'$ if $f(l)$ is lexicographically smaller than $f(l')$.
As a notation, we write $\perm{f_1\ f_2\ \dotsc}$ for the
permutation associated to the function $f=(f_1,f_2,\dotsc)$.
For example, $\rev \perm{f_1\ f_2\ \dotsc}=\perm{{-}f_1\ {-}f_2\ \dotsc}$.

As another example, if $X=Y=Z=[n]$ the four permutations
\begin{align*}
  \pi_1&=\perm{\phantom{-}x\ \ \phantom{-}y\ \ \phantom{-}z},\\
  \pi_2&=\perm{        {-}x\ \         {-}y\ \ \phantom{-}z},\\
  \pi_3&=\perm{        {-}x\ \ \phantom{-}y\ \         {-}z},\\
  \pi_4&=\perm{\phantom{-}x\ \         {-}y\ \         {-}z}.
\end{align*}
satisfy $\LCS(\pi_i,\pi_j)\leq n$ for $i\neq j$. Indeed, in any 
sequence $(x_1,y_1,z_1),(x_2,y_2,z_2),\dotsc$ that is increasing in 
both $\pi_i$ and $\pi_j$ only one of the coordinates can change. (This example
is from \cite{bh}).

\begin{proof}[Proof of \cref{thm:tweak}]
For $x\in \Z/p\Z$ let $\overline{x}$ be the element of $\{0,1,\dotsc,p-1\}$
that is congruent to $x$.
Let $X=Y=Z=\Z/p\Z$, and for each $i\in \Z/p\Z$ define the permutation 
\[
  \pi_i=\Bigl\lbrac\, \oline{i^2 x+iy+z}\quad \oline{2ix+y}\quad \oline{x}\,\Bigr\rbrac.
\]
We claim that $\LCS(\pi_i,\pi_j)\leq 4p-2$ for all $i\neq j$. 
Indeed, suppose $w$ is a common subsequence of $\pi_i$ and $\pi_j$.
Say, $w$ is the sequence $(x_1,y_1,z_1),(x_2,y_2,z_2),\dotsc \in (\Z/p\Z)^3$.
For $(I,J)\in \{0,1,\dotsc,p-1\}^2$ put
\[
  B(I,J)\eqdef \bigl\{(x,y,z)\in (\Z/p\Z)^3 : \overline{i^2 x+iy+z}=I,\ \overline{j^2 x+jy+z}=J\bigr\}.
\]
If $(x,y,z)\in B(I,J)$ we say that $(x,y,z)$ is in the \emph{bin} $(I,J)$.
Let $(I_1,J_1),(I_2,J_2),\dotsc$ be the sequence of bins 
for $(x_1,y_1,z_1),(x_2,y_2,z_2),\dotsc$. It is clear
that $I_1\leq I_2\leq \dotsb$ and $J_1\leq J_2\leq \dotsb$, and so at most
$2p-1$ bins are occupied. Hence, it suffices to prove that no more than
$2$ letters of $w$ fall into a same bin.

Since $i\neq j$, the two defining equations of $B(I,J)$ are linearly independent,
and so the set $B(I,J)$ is a line in $(\Z/p\Z)^3$. A computation yields that,
for fixed $(I,J)$, the set $B(I,J)$ is of the form $\bigl\{(x_0,y_0,z_0)+t(1,-i-j,ij) : t\in \Z/p\Z\bigr\}$ 
for some $(x_0,y_0,z_0)$. Thus, the set
\[
  \Bigl\{ \bigl(2ix+y,2jx+y\bigr) : (x,y,z)\in B(I,J) \Bigr\}
\]
is the line $\Bigl\{\bigl(2ix_0+y_0+(i-j)t,2jx_0+y_0+(j-i)t\bigr):t\in\Z/p\Z\Bigr\}$. As this line has slope $-1$,
the set
\[
  \Bigl\{ \bigl(\oline{2ix+y},\oline{2jx+y}\bigr) : (x,y,z)\in B(I,J) \Bigr\}
\]
is a union of at most two line segments of slope $-1$. Since
the sequence $w$ is increasing in both $\pi_i$ and $\pi_j$,
it follows that $w$ can contain at most two points from $B(I,J)$. Hence,
$\LCS_2(p,\P_{p^3})\leq 4p-2$. The inequality \eqref{eq:tweakone} 
follows by choosing $p$ to be the smallest prime exceeding $k^{1/3}$,
and noting that $p\leq k^{1/3}+O(k^{7/40})$ by \cite{prime_gaps}.

To derive the inequality \eqref{eq:tweakone} we start with a family of $t$ permutations
over a $t^3$-letter alphabet such that $\LCS$ of any two is at most
$4t+O(t^{21/40})$. We then select some $k$ letters of the alphabet,
and delete all the remaining letters from all the permutations.
\end{proof}

We extend the $\lbrac\dotsb\rbrac$ notation to the multipermutations.
The alphabet is still $X\times Y\times \dotsb $ for some sets
$X,Y,\dotsc$. There is now an additional set $R$,
and to each injective function $f\colon X\times Y\times\dotsb
\times R\to \Z\times\dotsb\times \Z$ we associate
a multipermutation in which each letter occurs $\abs{R}$
times, as follows. The set $R$ indexes copies
of the same letter in $\pi_f$. The occurrence of letter $l$ indexed by
$r\in R$ precedes the occurrence of letter
$l'$ indexed by $r'$ if $f(l,r)$ is lexicographically
smaller than $f(l',r')$.

\begin{proof}[Proof of \cref{thm:es_sharp}]
Let $X=[k_1]$, $Y=[k_2]$ and $R=[s]$. Consider the following two multipermutations:
\begin{align*}
\pi=\lbrac x\quad y\quad  \phantom{-}r\rbrac,\\
\pi'=\lbrac x\quad r\quad {-}y\rbrac.
\end{align*}
We first bound the $\LCS(\pi,\pi')$. Suppose
$(x_1,y_1,r_1),(x_2,y_2,r_2),\dotsc$ and 
$(x_1,y_1,r_1'),(x_2,y_2,r_2'),\dotsc$ are 
subsequences of $\pi$ and $\pi'$, respectively, that are equal 
as words. Then, for each $i$, either $x_{i+1}>x_i$, or
$x_{i+1}=x_i$ and $r_{i+1}>r_i$. Hence $\LCS(\pi,\pi')\leq k_1 s$.

We next bound $\LCS(\pi,\rev \pi')$. Consider a pair of sequences
as above. For each $i$ we must have $x_{i+1}=x_i$, and also
$y_{i+1}\geq y_i$ and $r_{i+1}'\leq r_i'$ with at least one of the inequalities
being strict. Hence $\LCS(\pi,\rev \pi')\leq k_2+s-1$.

Given any $n=ks$, let $k_1$ be the closest integer 
to $\sqrt{k/s}+\tfrac{1}{2}$, and let $k_2=\lceil k/k_1\rceil$. One can
then verify that $k_1 s\leq \sqrt{ks}+s$ and $k_2+s\leq \lceil \sqrt{ks}\rceil+s\leq \sqrt{ks}+s+1$.
\end{proof}

\begin{proof}[Proof of \cref{thm:multiperm_sharp}]
With hindsight let $X=[k_1]$, $Y=[k_2]$ and $Z=[k_3]$, where
$k_2$ is the closest integer to $(k/4s)^{1/3}+\tfrac{1}{3}$, and $k_1=2k_2$,
and $k_3=\lceil k/k_1k_2\rceil$. Put $R=[s]$. Consider the following four multipermutations:
\begin{align*}
\pi_1=\permt{+x \& +y \& +z \& +r},\\
\pi_2=\permt{-x \& -y \& +r \& +z},\\
\pi_3=\permt{-x \& +y \& -z \& +r},\\
\pi_4=\permt{+x \& -y \& +r \& -z}.
\end{align*}

To speed up the computation of $\LCS$ we observe that
$\LCS\left(\begin{smallmatrix}\perm{x\ f},\\\perm{x\ g}\phantom{,}\end{smallmatrix}\right)=\abs{X}\LCS\left(\begin{smallmatrix}\perm{f},\\\perm{g}\phantom{,}\end{smallmatrix}\right)$ and that
$\LCS\left(\begin{smallmatrix}\perm{\phantom{-}x\ f},\\\perm{-x\ g}\phantom{,}\end{smallmatrix}\right)=\LCS\left(\begin{smallmatrix}\perm{f},\\\perm{g}\phantom{,}\end{smallmatrix}\right)$.  
We thus have
\begin{align*}
\LCS(\pi_1,\pi_2)&=
  \LCS\left(\begin{smallmatrix}\perm{z\ r},\\\perm{r\ z}\phantom{,}\end{smallmatrix}\right)=\abs{Z}+s-1,\\
\LCS(\pi_1,\pi_4)&=
  \abs{X}\LCS\left(\begin{smallmatrix}\perm{z\ \phantom{-}r},\\\perm{r\ {-}z}\phantom{,}\end{smallmatrix}\right)=\abs{X}s,\\
 \LCS(\pi_1,\pi_3)&=\abs{Y}\LCS(\perm{r},\perm{r})=\abs{Y}s,\\
 \LCS(\pi_2,\pi_4)&=
  \abs{Y}\LCS\left(\begin{smallmatrix}\perm{r\ \phantom{-}z},\\\perm{r\ {-}z}\phantom{,}\end{smallmatrix}\right)=\abs{Y}(2s-1).
\end{align*}
The derivation of the equality $\LCS(\pi_3,\pi_4)=\abs{Z}+s-1$ is same as of $\LCS(\pi_1,\pi_2)$, whereas
the proof of $\LCS(\pi_2,\pi_3)=\abs{X}s$ is analogous to that of $\LCS(\pi_1,\pi_4)$.

With the choice of $k_1,k_2,k_3$ made above we have 
$(2s-1)\abs{Y}\leq \abs{X}s\leq (2s^2k)^{1/3}+\tfrac{5}{3}s$,
and $\abs{Z}+s-1\leq k/k_1k_2+s\leq k/2\bigl((k/4s)^{1/3}-1/6\bigr)^2=
(2s^2k)^{1/3}+\tfrac{5}{3}s+\frac{(1458k)^{1/3}-2s^{1/3}}{((54k)^{1/3}-s^{1/3})^2}\cdot \frac{s^{4/3}}{3}$,
which is at most $(2s^2k)^{1/3}+\tfrac{5}{3}s+s^{4/3}k^{-1/3}/3$ for $s\leq \tfrac{1}{5}k$.
\end{proof}

\section{Lower bounds}\label{sec:lowerbounds}
In this section we prove \cref{thm:lowerbound_perm} which gives the lower bound
on $\LCS_2(t,\P_k)$ for large $t$, and \cref{thm:multiperm_good} that shows
how to take advantage of repetitions in multipermutations. 

\begin{proof}[Proof of \cref{thm:lowerbound_perm}]\label{pf:lowerbound_perm}
Suppose that $\{\pi_1,\dotsc,\pi_t\}$ is any set of $t$ permutations, and 
put $L_{i,j}\eqdef\LCS(\pi_i,\pi_j)$.
For a letter $l\in [k]$ and $\pi\in \P_k$ let $\pi\{l\}$ be the 
prefix of $\pi$ that ends with $l$. For $i<j$ define 
a function $f_{i,j}\colon [k]\to [L_{i,j}]$ by $f_{i,j}(l)=\LCS\bigl(\pi_i\{l\},\pi_j\{l\}\bigr)$.

We say that a pair of letters $\{l,l'\}$
is \emph{nice} to the triple $i_1<i_2<i_3$ if
the three differences $f_{i_1,i_2}(l)-f_{i_1,i_2}(l')$,
$f_{i_1,i_3}(l)-f_{i_1,i_3}(l')$, and $f_{i_2,i_3}(l)-f_{i_2,i_3}(l')$
are either all negative, or all non-negative.\medskip

\begin{observation}
If $l$ and $l'$ are distinct, then there are
at least $\bigl(\tfrac{7}{16}-\tfrac{2}{t}\bigr)\binom{t}{3}$ triples $i_1<i_2<i_3$ for which $\{l,l'\}$
are nice.
\end{observation}
\begin{proof}[Proof of the observation]
Consider the complete graph on the vertex set $[t]$. Color its edge $i<\nobreak j$ red if $f_{i,j}(l)<f_{i,j}(l')$
and blue if $f_{i,j}(l)\geq f_{i,j}(l')$. A triple
$i_1<i_2<i_3$ is nice to $\{l,l'\}$ if $i_1i_2i_3$ is a monochromatic
triangle in the coloring.
Let $\mathcal{A}=\{i : l\text{ precedes }l'\text{ in }\pi_i\}$ and
$\mathcal{B}=\{i : l'\text{ precedes }l\text{ in }\pi_i\}$. Note that
all the edges in $\mathcal{A}$ are red, whereas all the 
edges in $\mathcal{B}$ are blue.

For an $i\in\mathcal{A}$ let $d(i)$ be the number of 
blue edges connecting $i$ to $\mathcal{B}$. For a $j\in\mathcal{B}$
let $d(j)$ be the number of red edges connecting $j$ to $\mathcal{A}$.
Finally let $M$ be the number of monochromatic triangles in our complete
graph. 
Since $\abs{\mathcal{A}}\abs{\mathcal{B}}=\sum_{i\in\mathcal{A}} d(i)+\sum_{j\in\mathcal{B}} d(j)
$ and $x\mapsto \binom{x}{2}$ is convex, it follows that
\begin{align*}
  M&=\binom{\abs{\mathcal{A}}}{3}+\binom{\abs{\mathcal{B}}}{3}+\sum_{i\in\mathcal{A}} \binom{d(i)}{2}+\sum_{j\in\mathcal{B}} \binom{d(j)}{2}.\\
   &\geq \binom{\abs{\mathcal{A}}}{3}+\binom{\abs{\mathcal{B}}}{3}+
    t\binom{\abs{\mathcal{A}}\abs{\mathcal{B}}/t}{2}.
\end{align*}
Let $x=\abs{\mathcal{A}}$. Then $\abs{\mathcal{B}}=t-x$. With a bit of calculus
we can compute the derivative of the right-hand side of the inequality 
above with respect to $x$ to be $(x/t-1/2)\bigl(t(t-1)-2xt+2x^2\bigr)$,
from which it follows that the minimum is at $x=n/2$, and so
\[
  M \geq 2\binom{t/2}{3}+t\binom{t/4}{2}\geq \left(\frac{7}{16}-\frac{2}{t}\right)\binom{t}{3}.\qedhere
\]
\end{proof}

The observation implies that there is a triple $i_1<i_2<i_3$ that is nice to $(\frac{7}{16}-\frac{2}{t})\binom{k}{2}$
pairs $\{l,l'\}\in \binom{[k]}{2}$. Note, however, that the number of unordered pairs $\{(x,y,z),(x',y',z')\}$
such that $x\leq x',y\leq y',z\leq z'$ in a box $[L_{i,j}]\times [L_{i,k}]\times [L_{j,k}]$ is less than
$
  \binom{L_{i_1,i_2}+1}{2}\binom{L_{i_1,i_3}+1}{2}\binom{L_{i_2,i_3}+1}{2}
$.
Hence,
\[
  \left(\frac{7}{16}-\frac{2}{t}\right)\binom{k}{2}\leq \binom{L_{i_1,i_2}+1}{2}\binom{L_{i_1,i_3}+1}{2}\binom{L_{i_2,i_3}+1}{2}.
\]
Since $(k-1)^2/2\leq \binom{k}{2}$ and $\binom{L+1}{2}\leq (L+1)^2/2$, we have 
$\max(L_{i_1,i_2},L_{i_1,i_3},L_{i_2,i_3})\geq \left(\frac{7}{4}-\frac{8}{t}\right)^{1/6}k^{1/3}-\nobreak 2$.
\end{proof}

For a proof of \cref{thm:multiperm_good} we need a lemma about monotone functions 
that is of independent interest.
A function $f(x,y)$ of two variables is said to be \emph{strongly monotone} if
the inequalities $f(x,y)\leq f(x',y)$, $f(x,y)\leq f(x,y')$ and $f(x,y)<f(x',y')$ hold 
whenever $x<x'$ and $y<y'$.
\begin{lemma}\label{lem:strongly_monotone}
Suppose $f_1,f_2,f_3\colon [s]^2\to\Z$ are strongly monotone functions,
and 
\[
f(x,y,z)=\bigl(f_1(x,y),f_2(x,z),f_3(y,z)\bigr).
\] Then
$f$ takes at least $s^2/6$ distinct values on $[s]^3$.
\end{lemma}
\begin{proof}
The key observation is that if $f(x,y,z)=f(x',y',z')$, then
$(x,y,z)$ and $(x',y',z')$ agree in at least one coordinate.
Indeed, if it were not true, then by swapping $p$ with $p'$ and 
renaming the coordinates if necessary, we could have arranged
that $x<x'$ and $y<y'$, which would have contradicted the strong monotonicity of $f_1$.

Let $(x_0,y_0,z_0)\in [s]^3$ be arbitrary, and consider the 
set $E\eqdef \{(x,y,z) : f(x,y,z)=f(x_0,y_0,z_0)\}$.
The observation tells us that the set $E$ is contained in a union of 
three coordinate hyperplanes, namely 
$H_x=\{ (x_0,y,z) : y,z\in [s]^2\}$, 
$H_y=\{ (x,y_0,z) : x,z\in [s]^2\}$,
and
$H_z=\{ (x,y,z_0) : x,y\in [s]^2\}$.
We claim that each of these hyperplanes contains at 
most $2s-1$ points of $E$. Indeed,
$H_x\cap E$ cannot contain
two points $(x_0,y,z)$ and $(x_0,y',z')$ such that
$y<y'$ and $z<z'$, for that would have contradicted strong monotonicity of 
the function $f_3$. In particular $y-z$ is distinct
as $(x_0,y,z)$ runs over points of $H_x\cap E$.
Thus, $\abs{E}\leq 3(2s-1)\leq 6s$. Since $(x_0,y_0,z_0)$
was arbitrary, the image of $f$ must be of size at least $s^3/6s$,
which completes the proof.
\end{proof}

\begin{proof}[Proof of \cref{thm:multiperm_good}]
Put $L_{i,j}\eqdef \LCS(\pi_i,\pi_j)$.
Let $\mathcal{I}$ be the set of all pairs $(l,p)$ consisting of a letter $l\in [k]$
and an integer $1\leq p\leq s_l$. For $(l,p)\in \mathcal{I}$
denote by $\pi_i \{ l, p \}$ the prefix of the multipermutation $\pi_i$ that ends with the $p$'th copy of the letter $l$.  Put $\mathcal{I}_3\eqdef \bigl\{ (l,p_1,p_2,p_3) : (l,p_i)\in \mathcal{I} \bigr\}$.
Define a function 
$f\colon \mathcal{I}_3\to [L_{1,2}]\times [L_{1,3}]\times [L_{2,3}]$ by
\[
  f(l,p_1,p_2,p_3)=\bigl(\LCS(\pi_1 \{ l, p_1 \},\pi_2 \{ l, p_2 \}),\,\LCS(\pi_1 \{ l, p_1 \},\pi_3 \{ l, p_3 \}),\,\LCS(\pi_2 \{ l, p_2 \},\pi_3 \{ l, p_3 \})\bigr).
\]
If $l,l'\in [k]$ are two different letters, then 
$f(l,p)\neq f(l',p')$ for all $p\in [s_l]^3$ and $p'\in [s_{l'}]^3$. 
Indeed, interchanging the roles
of $(l,p)$ and $(l',p')$ and renaming the multipermutations if needed, we may assume that
$\pi_1\{l,p_1\}$ and $\pi_2\{l,p_2\}$ are longer than
$\pi_1\{l',p_1'\}$ and $\pi_2\{l',p_2'\}$ respectively, and so $f_1(l,p)>\nobreak f_1(l',p')$
because the longest common subsequence between $\pi_1\{l',p_1'\}$ and $\pi_2\{l',p_2'\}$
can be extended to a longer common subsequence of $\pi_1\{l,p_1\}$ and $\pi_2\{l,p_2\}$.

Hence, $f(l,[s_l]^3)\cap f(l',[s_l']^3)=\emptyset$ for distinct $l,l'$. However,
the preceding lemma shows that for any fixed letter $l$ 
we have $\abs{f(l,[s_l]^3)}\geq s_l^2/6$. Thus $\tfrac{1}{6}\sum s_l^2\leq L_{1,2}L_{1,3}L_{2,3}$
and the theorem follows.
\end{proof}

\section{Concluding remarks}
\begin{itemize}
\item The upper bounds on $\LT(k,n)$ both in this paper and in \cite{app_twins} come from
random words. Estimating $\LT(w)$ for a random $w\in [k]^n$ is an interesting problem
on its own. A result of Kiwi, Loebl and Matou\v{s}ek \cite{kiwi_loebl_matousek}
asserts that if $w_1,w_2$ are two random words of length $n/2$, then
$\E[\LCS(w_1,w_2)]\sim \frac{1}{\sqrt{k}}n$ for large $k$. Hence
$1\leq \lim_n \frac{\E[\LT(w)]}{n\sqrt{k}}\leq e$ for large $k$.

\item In \cite{app_twins} a more general problem has been considered: Instead of looking for twins in $w$,
one can seek $T$-tuplets, which are $T$-tuples of disjoint subsequences that are equal as words. One
then defines the quantity $\LT_T(w)$ to be the length of longest $T$-tuplets in $w$, and defines
$\LT_T(k,n)$ in a manner analogous to $\LT(k,n)$.
\Cref{thm:klarge} extends easily to $\LT_T(k,n)\geq C_T k^{-1+1/\binom{2T-1}{T}}n$
with the help of the following estimate:
\begin{theorem}
For every $T$ and $k$ we have $\LCS_T(2T-1,\P_k)\geq k^{1/\binom{2T-1}{T}}$.
\end{theorem}
\begin{proof}
For a letter $l\in [k]$ and any set of permutations $P\subset \P_k$
we denote by $P\{l\}$ the prefixes of all $\pi\in P$ that end with
the letter $l$.  Suppose $\Pi\subset \P_k$ is any set of $2T-1$ permutations.
For each $I\subset \binom{\Pi}{T}$ and $l\in [k]$ put
$f_I(l)\eqdef \LCS(I\{l\})$. Define the function $f\colon [k]\to \mathbb{N}^{\binom{\Pi}{T}}$
by $f(l)_I=f_I(l)$. If $l,l'\in [k]$ are distinct letters, then by the pigeonhole principle 
there are at least $T$ permutations in which the relative order of $l$ and $l'$ is the same,
and so $f(l)\neq f(l')$. So, $f$ is injective, and thus its image cannot be
contained in a box with of length smaller than\nobreak\ $k^{1/\binom{2T-1}{T}}$.
\end{proof}

The preceding theorem is sharp. First, $2T-1$ cannot be reduced to $2T-2$
because of a family that consists of $T-1$ copies of $\wordf{123}\dotsb \wordf{k}$ and $T-1$
copies of $\wordf{k}\dotsb \wordf{321}$. Second, the bound itself cannot be improved, for
$\LCS_T(2T,\P_k)\leq k^{1/\binom{2T-1}{T}}$ as the following example shows:
For each $I\subset \binom{[2T-1]}{T}$ associate a set $X_I=[k]$, and 
define $f_{i,I}\colon X_I\to \Z$ by $f_{i,I}(x_I)=x_I$ if $i\in I$
and $f_{i,I}(x_I)=-x_I$ if $i\not\in I$. Also define $f_{0,I}$ by $f_{0,I}(x_I)=x_I$.
Define $f_i=(f_{i,I})_I$, which is a function from $\prod_I X_I$ to $\Z^{\binom{2T-1}{T}}$,
and put $\pi_i=\lbrac f_i \rbrac$ as in \cref{sec:constr}.
It is then easy to check that $\Pi=\{\pi_0,\pi_1,\dotsc,\pi_{2T-1}\}$ satisfies 
$\LCS_T(\Pi)=k$.

\item In this paper we neglected to address the most basic extremal problem
on $\LCS$ in sets of words, the estimation of $\LCS_2(t,[k]^n)$. Denote by $l^n$ the
word made of $n$ copies of the letter $l$. The family
$\{\wordf{1}^n,\wordf{2}^n,\dotsc,\wordf{k}^n\}$ shows that $\LCS_2(k,[k]^n)=0$. By
labeling a word with the most popular letter that occurs in it, and applying
the pigeonhole principle, one derives $\LCS_2(k+1,[k]^n)\geq n/k$. The bound is sharp, 
as the family of $k+2$ words
$\{\wordf{1}^n,\wordf{2}^n,\dotsc,\wordf{k}^n,
\wordf{1}^{n/k}\wordf{2}^{n/k}\dotsb\wordf{k}^{n/k},
\wordf{k}^{n/k}\dotsb\wordf{2}^{n/k}\wordf{1}^{n/k}
\}$ shows. For fixed $t\geq k+3$ it is possible to prove that
$\LCS_2(t,[k]^n)\leq \bigl(1+o(1)\bigr)n/k$ by considering
a family of words of the form 
$w_m=\bigl(\wordf{1}^m\dotsb\wordf{k}^m\bigr)^{n/mk}$
for a quickly growing sequence of values of $m$.
Indeed, if $m_1<m_2$, then a common subsequence of $w_{m_1}$ and $w_{m_2}$
is necessarily of the form $l_1^{p_1}l_2^{p_2}\dots l_r^{p_r}$
for some $r\leq n/m_2$, but each subsequence of form
$l^p$ in $w_1$ spans an interval of length at least $\lfloor \frac{p-1}{m_1}\rfloor m_1 k\geq (p-m_1)k$
in $w_1$. Since $\sum (p_i-m_1)k\leq n$, we obtain that $\LCS_2(w_{m_1},w_{m_2})\leq \bigl(\frac{1}{k}+\frac{m_1}{m_2}\bigr)n$.
In a recent work of the first author with Jie Ma \cite{bukh_ma}, we showed that this construction is essentially optimal 
when the number of words is constant.

\item It is hard to resist the conjecture that the correct bound in \cref{lem:strongly_monotone}
is $s^2$ in place of $s^2/6$. If true, it would be sharp in view of the function $f(x,y,z)=(x,x,y)$.
\end{itemize}

\paragraph{Acknowledgments.} The first author gratefully acknowledges helpful discussions with
Yury Person and Sevak Mkrtchyan. We thank the referees for suggestions that helped to improve
the exposition. All the errors remain ours.

\bibliographystyle{plain}
\bibliography{llcs}

\end{document}